\documentclass{elsart}

\usepackage{graphicx}
\usepackage{amssymb}
\usepackage[latin1]{inputenc} 

\begin{document}

\begin{frontmatter}

% Title, authors and addresses

% use the thanksref command within \title, \author or \address for footnotes;
% use the corauthref command within \author for corresponding author footnotes;
% use the ead command for the email address,
% and the form \ead[url] for the home page:
% \title{Title\thanksref{label1}}
% \thanks[label1]{}
% \author{Name\corauthref{cor1}\thanksref{label2}}
% \ead{email address}
% \ead[url]{home page}
% \thanks[label2]{}
% \corauth[cor1]{}
% \address{Address\thanksref{label3}}
% \thanks[label3]{}

\title{Perfect Skolem sets}

% use optional labels to link authors explicitly to addresses:
% \author[label1,label2]{}
% \address[label1]{}
% \address[label2]{}

\author{Gustav Nordh \thanksref{cugs}}

\address{Department of Computer and Information Science, Linköpings Universitet, \\S-581 83 Linköping, Sweden}
\ead{gusno@ida.liu.se}
\thanks[cugs]{Supported by the {\em National Graduate School in Computer Science} (CUGS), Sweden.}

\begin{abstract}
% Text of abstract
A Skolem sequence is a sequence $a_1,a_2,\dots,a_{2n}$ (where $a_i \in A = \{1 \dots n \})$, each $a_i$ 
occurs exactly twice in the sequence and the two occurrences are exactly $a_i$ positions apart. A set $A$ that 
can be used to construct Skolem sequences is called a Skolem set. The problem 
of deciding which sets of the form $A = \{1 \dots n\}$ are Skolem sets was solved
by Thoralf Skolem in the late 1950's. We study the natural generalization where $A$ is allowed to be any set of $n$ positive integers. We give necessary conditions for the existence of Skolem sets of this generalized form. We conjecture these necessary conditions to be sufficient, and give computational evidence in favor of our conjecture.
We investigate special cases of the conjecture and prove that the conjecture hold for some of them. We also study enumerative questions
and show that this problem has strong connections with problems related to permutation displacements.
\end{abstract}

\end{frontmatter}

% main text
%\section{}
%\label{}

% The Appendices part is started with the command \appendix;
% appendix sections are then done as normal sections
% \appendix

% \section{}
% \label{}

%\documentclass[twoside]{article}
%\usepackage[latin1]{inputenc}   % Enter your text in ISO-Latin 1
%\usepackage[dvips]{graphicx}
%\usepackage{latexsym}
%\author{Gustav Nordh\\Department of Computer and Information Science\\Linköpings Universitet\\S-581 83 Linköping, Sweden\\
%Email: gusno@ida.liu.se}
%\title{Perfect Skolem sequences}
%\pagestyle{empty}
%\date{6 Mars, 2003}

%\begin{document}
%\maketitle
%
%\vspace{7mm}
%\centerline{\bf Abstract}
%{\small A Skolem sequence is a sequence $a_1,a_2,\dots,a_{2n}$ (where $a_i \in A = \{1 \dots n \})$, each $a_i$ 
%occurs exactly twice in the sequence and the two occurrences are exactly $a_i$ positions apart. A set $A$ that 
%can be used to construct Skolem sequences is called a Skolem set. The problem 
%of deciding which sets of the form $A = \{1 \dots n\}$ that are Skolem sets was solved
%by Th. Skolem in the late 1950's. We study the natural generalization where $A$ is allowed to be any set (or even multiset) of $n$ positive integers. We give necessary conditions for the existence of Skolem sets of this generalized form. We conjecture these necessary conditions to be sufficient, and give computational evidence in favor of our conjecture.
%We also look at enumerative aspects and prove that the number of multi Skolem sequences (multisets allowed) of order $n$ is $(2n-1)!!$.}
%\vskip 10pt
%\textbf{Keywords:} Skolem sequence, design theory, Langford's problem
%\newpage
\section{Introduction}
{\em Skolem sequences\/} were introduced by Thoralf Skolem in 1957 \cite{Skol}, for the construction of Steiner triple systems. He looked at sets of the
form $A=\{1,2,\dots,n\}$ and asked whether you always could form a sequence with two copies of every element $k$ in 
the set so that the two copies
of $k$ were placed $k$ places apart in the sequence. For example, the set $\{1,2,3,4\}$ can be used to form the sequence 
$42324311$, but the set
$\{1,2,3\}$ cannot be used to form such a sequence. Close relatives to Skolem sequences are Langford sequences. 
%Langford sequences
%are equivalent to Skolem sequences with the difference that the two copies of $k$ must be placed $k+1$ places apart. For example the
%set $\{1,2,3\}$ can be used to construct the sequence $312132$. Dudley Langford published an article in 1958 \cite{Lang} where he %managed to
%construct such sequences from the sets $\{1 \dots n\}$ with $n = 5,6,9$ and $10$. 
A {\em Langford sequence\/}
is a Skolem sequence which starts late. That is, a sequence that is constructed from a set of the form $\{a,a+1, \dots, b-1,b\}$.
Dudley Langford published an article in 1958 \cite{Lang} where he managed to
construct such sequences from the sets $\{2 \dots n\}$ with $n = 5,6,9$ and $10$.

Since these early papers appeared many different aspects of Skolem sequences and Langford sequences have been studied. One reason for
them being so well studied is that they have important applications in several branches of mathematics, most notably
in design theory and graph labeling; see \cite{HCD}. 
%One real world application of the generalization of Skolem sequences we are about to study is the generation of binary sequences with %controllable complexity \cite{Gorth}.   

A natural generalization of Skolem sequences that is not very well studied is when we allow
the set of integers used to generate the sequences to be any set or multiset of positive integers. We call such sequences
{\em perfect Skolem-type sequences}. This generalization is the main topic of this article.    

We call a set that can be used to construct a perfect Skolem-type sequence a {\em perfect Skolem set}. For example the set
$\{2,3,5,6\}$ is a perfect Skolem set since it can be used to generate the sequence $56232536$. When we want to point out that the sequences are constructed from multisets we call the
sequences {\em perfect multi Skolem-type sequences\/} and the corresponding multiset is a {\em perfect multi Skolem set}.

The problem of deciding whether the set $P=\{1,2,\dots,2n\}$ can be partitioned into the differences 
in $A=\{a_1,a_2,\dots,a_n\}$ is exactly the problem of deciding whether $A$ is a perfect multi Skolem set.
We can generalize the notion of Skolem sequences even further by allowing the set of positions $P$ to be an arbitrary set.
If the set of positions $P$ can be partitioned into the differences $A$, then we say that $(P,A)$ is {\em generalized multi Skolem\/} 
(if $A$ is not allowed to be a multi set, then we say that $(P,A)$ is {\em generalized Skolem\/}).
For example, $P=\{1,2,4,5,7,8\}$, $A=\{1,6,6\}$ is generalized multi Skolem by the sequence $66\_11\_66$.
Given a set $P$ and a multiset $A$, the problem of deciding whether $(P,A)$ is generalized multi Skolem is NP-complete \cite{Nordh}.

We now state some of the most important known results that are related to the generalization that we are about to study.  

\begin{thm}[\cite{Skol}]
$\{1,2,3,4,\dots,n-1,n\}$ is a perfect Skolem set if and only if $n \equiv 0,1 \pmod{4}$.
\label{skolem}
\end{thm}
The preceding theorem classify those sets which can be used to construct Skolem sequences.

\begin{thm}[\cite{BBG76,Simp}]
A set $A = \{a,a+1,a+2,\dots,b-1,b\}$ is a perfect Skolem set if and only if
\begin{enumerate}
\item $|A| \geq 2a -1$; and
\item $|A| \equiv 0,1 \pmod{4}$ when $a$ is odd, $|A| \equiv 0,3 \pmod{4}$ when $a$ is even.
\end{enumerate}
\label{langford}
\end{thm}
The preceding theorem classify those sets which can be used to construct Langford sequences.

\begin{thm}[\cite{Shal}]
A set $\{1,2,\dots,m-1,m+1,\dots,n\}$ is a perfect Skolem set if and only if $n\equiv 0,1 \pmod{4}$ and $m$ is odd, 
or $n \equiv 2,3 \pmod{4}$ and $m$ is even.
\label{nearskolem}
\end{thm}
Sequences constructed from these sets in the preceding theorem are known in the literature as near-Skolem sequences.

\begin{thm}[\cite{Bak94}]
A multiset $\{1^m,2^m,\dots,n^m\}$ is a perfect multi Skolem set if and only if
\begin{enumerate}
\item $n \equiv 0,1 \pmod{4}$; or
\item $n \equiv 2,3 \pmod{4}$ and $m$ is even.
\end{enumerate}
\label{mfold}
\end{thm}
Note that by $a^b$ we mean $b$ copies of $a$.
Sequences constructed from these sets in the preceding theorem are known in the literature as $m$-fold Skolem sequences.

\begin{thm}[\cite{Bak}]
$P = \{1,2,\dots,k-1,k+1,k+2,\dots,2n+1\}$ and $A=\{1,2,\dots,n\}$ is generalized Skolem if and only if $n \equiv 0,1 \pmod{4}$ for
$k$ odd and $n \equiv 2,3 \pmod{4}$ for $k$ even.
\label{kext}
\end{thm}
Sequences constructed from these sets in the preceding theorem are known in the literature as $k$-extended Skolem sequences.

\section{Necessary conditions and conjectures}
We give necessary conditions for a set $A$ to be a perfect Skolem set. These necessary conditions are then conjectured to be 
sufficient for a set $A$ to be a perfect Skolem set. We first state the necessary conditions in a more general form.

\begin{thm}
Let $A = \{a_1,a_2,\dots,a_n\}$ be a multiset of differences with $a_1 \geq a_2 \geq \dots \geq a_n$, and let 
$P = \{p_1,p_2,\dots,p_{2n}\}$, with $p_1 < p_2 < \dots < p_{2n}$. If $P$ can be partitioned into the differences in $A$, then

\begin{equation}
\label{eqparity}
\sum^n_{i=1} a_i \equiv \sum^{2n}_{i=1} p_i \;\; \pmod{2},
\end{equation}
\begin{equation}
\label{eqdensity}
\sum^m_{i=1} a_i \leq \left( \sum^{2n}_{i = 2n+1 -m} p_i \right) - \left( \sum^m_{i=1} p_i \right), \;\;\; (m = 1,\dots,n).
\end{equation}

\label{necgeneral}
\end{thm}
\begin{pf}
Let $a_i = s_i - t_i$, ($i = 1,2,\dots,n$). Then $a_i \equiv s_i + t_i \pmod{2}$ and summing over $i$ gives (\ref{eqparity}).
Also for each $1 \leq m \leq n$ we have 
$a_1 + \dots + a_m = (\sum^m_{i=1} s_i) - (\sum^m_{i=1} t_i) \leq (\sum^{2n}_{i=2n+1-m} p_i) - (\sum^m_{i=1} p_i)$, and 
(\ref{eqdensity}) is proved.
\qed
\end{pf}
Conditions \ref{eqparity} and \ref{eqdensity} are referred to as parity and density conditions respectively. Note that it is easy to see
that these conditions are not sufficient for $(P,A)$ to be generalized Skolem. For example $P=\{1,2,4,5\}$ and $A=\{1,3\}$ satisfies 
both \ref{eqparity} and \ref{eqdensity} but is not generalized Skolem.

We leave it up to the reader to verify the following consequence of the preceding theorem.
\begin{cor}
If $P = \{1,2,\dots,2n\}$, then the two necessary conditions in the preceding theorem 
reduce to (\ref{eqparity}) the number of even $a_i$'s is even, and (\ref{eqdensity}) $\sum^m_{i=1} a_i \leq m(2n-m)$ for each 
$1 \leq m \leq n$.
\label{necperf}
\end{cor}

Surprisingly, when $P=\{1,2,\dots,2n\}$ and $A$ is an ordinary set (i.e., not a multiset), then the necessary conditions 
seems to be sufficient.
\begin{conj}
A set $A = \{a_1,a_2,\dots,a_n\}$ with $a_1 > a_2 > \dots > a_n$ is a perfect Skolem set if and only if the number of even $a_i$'s is even, and $\sum^m_{i=1} a_i \leq m(2n-m)$ for each 
$1 \leq m \leq n$.
\label{perfconj}
\end{conj}
The conjecture has been verified by computer search for all sets of cardinality $20$ or less. The conjecture 
does not hold when we allow the 
sets to be multisets. A minimal counterexample is the multiset $\{1,3,3\}$, it passes both necessary conditions, but cannot be used 
to construct a perfect multi Skolem-type sequence. The existence results for Skolem sequences, Langford sequences, and near-Skolem sequences in Theorems \ref{skolem}, \ref{langford}, and \ref{nearskolem} are special cases where we know that the conjecture holds. Many new minor special
cases where the conjecture holds can be easily deduced from already known results on, e.g., $k$-extended Skolem sequences.
For example, the results in \cite{Bak} imply that $A=\{1,2,\dots,n\} \cup \{2j+1\}$ is a perfect Skolem set for $(n-1)/2<j \leq n$ if $n \equiv 0,1 \pmod{4}$, and that $A=\{1,2,\dots,n\} \cup \{2j\}$ is a perfect Skolem set for $n/2<j \leq n$ if $n \equiv 2,3 \pmod{4}$.
The results in \cite{Lin01} imply that $\{1,3,5,\dots,2n-1\} \cup \{2,2j\}$ is a perfect Skolem set for $n \geq 2$ and $2 \leq j \leq n$. We give more evidence for the conjecture in terms of special cases where it can be proved to hold in Section \ref{sectionconstr}.

It is easy to formulate interesting special cases of Conjecture \ref{perfconj}. Consider for example sets $A$ of the form $A \subset \{1,2,\dots,n\}$, 
where $|A| = n-2$, i.e., $A$ is a subset of $\{1,2,\dots,n\}$ such that exactly two elements $a_i,a_j \in \{1,2,\dots,n\}$ are 
missing from $A$. 
\begin{conj}
A subset $A$ of $\{1,2,\dots,n\}$ such that exactly two elements $a_i,a_j \in \{1,2,\dots,n\}$ are missing from $A$, is a perfect 
Skolem set if and only if the following conditions hold.
\begin{enumerate}
\item If $n \equiv 0$ or $1$ $(mod$ $n)$ then $a_i \equiv a_j$ $(mod$ $2)$.    
\item If $n \equiv 2$ or $3$ $(mod$ $n)$ then $a_i \equiv a_j-1$ $(mod$ $2)$.
\item $A \notin \{\{3\},\{2,4\},\{2,4,5\},\{3,4,5,6\}\}$. 
\end{enumerate}
\end{conj}
%Even a proof of this special case of Conjecture \ref{perfconj} would be interesting. 
%The conjecture has been verified by computer 
%search for all sets of size $13$ or less. 
Note that conditions $1$ and $2$ corresponds to the parity condition in Theorem \ref{necperf}, 
condition $3$ corresponds to the density condition in Theorem \ref{necperf}. The sets in condition $3$ are the only ones 
(of this particular 
form) that have the right parity but do not satisfy the density condition in Theorem \ref{necperf}.

Another particularly interesting special case of Conjecture \ref{perfconj} emerges when we add the condition 
$$\sum^n_{i=1} a_i = \left( \sum^{2n}_{i = 2n+1 -m} p_i \right) - \left( \sum^m_{i=1} p_i \right)$$
to the necessary conditions in Theorem \ref{necgeneral}. We call a multiset $A$ satisfying these conditions 
{\em extremal}, and the corresponding sequences for {\em extremal (multi) Skolem-type sequences}. 
These sets are extremal in the sense
that adding $1$ (or more) to any of the elements in $A$ would force $A$ to violate the density condition.
When $P=\{1,2,\dots,2n\}$ and $A$ is an ordinary set (i.e., not a multiset), then we end up with the following special case of Conjecture
\ref{perfconj}.
\begin{conj}
A set $A = \{a_1,a_2,\dots,a_n\}$ with $a_1 > a_2 > \dots > a_n$ and $\sum^n_{i=1} a_i = n^2$ is a perfect extremal Skolem set if and
only if the number of even $a_i$'s is even, and $\sum^m_{i=1} a_i \leq m(2n-m)$ for each 
$1 \leq m \leq n$.
\label{extconj}
\end{conj}
Note that again the conjecture does not hold when we allow $A$ to be a multiset. A minimal counter example is $\{4,4,4,8,8,8\}$. 
As we will see in Section \ref{restrictedpermutations}, extremal sets and sequences have interesting connections to problems related
to permutation displacements.

\section{Enumerative aspects}
The problem of deciding how many Langford sequences of a given order there are goes under the name Langford's problem. It is quite well studied, but so far no closed
formula is known and the fastest known algorithm runs in exponential time. For an overview of the results; see 
\cite{Miller}. Recently Langford's problem has began to be used as a benchmark by the computer science community to evaluate various Constraint Satisfaction techniques \cite{Habbas}.
The problem of deciding how many Skolem sequences of a given order there are has not been given the same attention. A table of the number of Skolem sequences of orders $1$ to $13$ is available in \cite{HCD}. We extend this table by computing the number of Skolem 
sequences of orders $16$ and $17$ (the actual numbers are presented in Table \ref{numberofskol}).

\begin{table}[!h]
%\begin{center}
\caption{Number of Skolem sequences.}
\label{numberofskol}
\begin{tabular*}{12cm}{|c@{\extracolsep{\fill}}|cccccc|}
\hline
Order & 1 & 4 & 5 & 8 & 9 & 12 \\
\hline
$\#$Sequences $\;$ & 1 & 6 & 10 & 504 & 2656 & 455936 \\
\hline
\end{tabular*}
\begin{tabular*}{12cm}{|c@{\extracolsep{\fill}}|ccc|}
\hline
Order & 13 & 16 & 17  \\
\hline
$\#$Sequences $\;$ & 3040560 & 1400156768 & 12248982496\\
\hline
\end{tabular*}
\end{table}

As for Langford's problem, no closed formula is known and the fastest known algorithm runs in exponential time.
We begin our investigation by giving a closed formula for the number of perfect multi Skolem-type sequences of any given order. 
Note that we say that
a sequence is of order $n$ if the set that it is generated from has cardinality $n$.
\begin{thm}
The number of perfect multi Skolem-type sequences of order $n$ is $(2n-1)!!$.
\label{nrmulti}
\end{thm}
\begin{pf}
Think of the sequence as a tape with $2n$ empty cells. Pick $2$ cells, say $i$ and $j$, from the $2n$ available cells. Write the number
$|i-j|$ in cells $i$ and $j$. Now repeat the procedure and choose $2$ cells from the remaining $2n-2$ cells. Continue in this way 
until no more empty cells remains. It is clear that the sequence of numbers on the tape is a perfect multi Skolem-type sequence and that 
every perfect multi Skolem-type sequence can be constructed in this way. Note that we are not interested in the order in which the $n$ 
pairs are chosen, we do not care for example if $2$ and $6$ were chosen before or after $1$ and $3$. We compensate for this by 
putting $n!$ in the denominator. Thus, the number of perfect multi Skolem-type sequences of order $n$ is
$$\frac{{2n \choose 2}{2n-2 \choose 2} \dots {4 \choose 2}{2 \choose 2}}{n!} = \frac{2n!}{n!2^n} = (2n-1)(2n-3) \dots (3)(1)=(2n-1)!!$$
\qed
\end{pf}

By using a similar argument as in the proof of the preceding theorem, we can also give a closed formula for the number of
perfect extremal multi Skolem-type sequences. An important observation about perfect extremal (multi) Skolem-type sequences is that
the first occurrence of each element $a_i$ in $A$ must be placed in the left half of the sequence and the second occurrence must be
placed in the right half of the sequence. This is because in any partition $(s_i, t_i)$, $i \in \{1,\dots,n\}$ of $\{1,2,\dots,2n\}$ into 
the differences in $\{s_i - t_i\} = A = \{a_1,a_2,\cdot, a_n\}$ satisfying $\sum^n_{i=1} a_i = n^2$, we must have $t_i \leq n$ and 
$s_i > n$ for all $i \in \{1,\dots,n\}$.
\begin{thm}
The number of perfect extremal multi Skolem-type sequences of order $n$ is $n!$.
\label{nrmaxmulti}
\end{thm}
\begin{pf}
%Note that perfect extremal multi Skolem-type sequences have the property that the first occurrence of each element $a_i$ in $A$ is %placed in the left half of the sequence and the second occurrence is placed in the right half of the sequence.
In analogy with the proof of the
preceding theorem, we pick one cell from the $n$ available cells in the right half of the sequence (tape), say $i$, and write the number
$i-1$ in cells $1$ and $i$. We repeat the procedure and choose one cell from the remaining $n-1$ cells in the right half 
of the sequence, say $j$, and write the number $j-2$ in cells $2$ and $j$. We continue in this way until no more empty cells remains.
It should be clear that the sequence of numbers on the tape is a perfect extremal multi Skolem-type sequence and that every perfect
extremal multi Skolem-type sequence can be constructed in this way. Thus, the number of perfect extremal multi Skolem-type sequences of 
order $n$ is $n!$.
\qed
\end{pf}

Despite our success in finding closed formulas for the number of perfect multi Skolem-type sequences and perfect 
extremal multi Skolem-type sequences, we have not been as 
fortunate when it comes to finding a 
formula for the number of perfect Skolem-type sequences (those constructed from ordinary sets). We have put considerable time and 
effort 
%(approximately 4 months of computation on a 350 MHz Pentium II computer) 
in computing the number of perfect Skolem 
sequences of order at most $13$. The actual numbers are available in Table \ref{nrskol} and Table \ref{nrextskol}. In Section 
\ref{restrictedpermutations} 
we give nice interpretations of these numbers in terms of permutation matrices. 

%\begin{table}[!h]
%%\begin{center}
%\caption{Number of perfect Skolem sequences.}
%\begin{tabular}{|c|l|}
%\hline
%Order & $\#$perfect Skolem sequences \\
%\hline
%1 & 1 \\
%2 & 1 \\
%3 & 5 \\
%4 & 29 \\
%5 & 145 \\
%6 & 957 \\
%7 & 8397 \\
%8 & 85169 \\
%9 & 944221 \\
%10 & 11639417 \\
%11 & 160699437 \\
%12 & 2430145085 \\
%13 & 39776366397 \\
%\hline
%\end{tabular}
%%\end{center}
%\label{nrskol}
%\end{table}

\begin{table}[!h]
%\begin{center}
\caption{Number of perfect Skolem-type sequences.}
\begin{tabular*}{13cm}{|c@{\extracolsep{\fill}}|ccccccccc|}
\hline
Order & 1 & 2 & 3 & 4 & 5 & 6 & 7 & 8 & 9 \\
\hline
$\#$Sequences $\;$ & 1 & 1 & 5 & 29 & 145 & 957 & 8397 & 85169 & 944221\\
\hline
\end{tabular*}
\begin{tabular*}{13cm}{|c@{\extracolsep{\fill}}|cccc|}
\hline
Order & 10 & 11 & 12 & 13  \\
\hline
$\#$Sequences $\;$ & 11639417 & 160699437 & 2430145085 & 39776366397 \\
\hline
\end{tabular*}
%\end{center}
\label{nrskol}
\end{table}

%\begin{table}[!h]
%%\begin{center}
%\caption{Number of perfect extremal Skolem sequences.}
%\begin{tabular}{|c|l|}
%\hline
%Order & $\#$perfect extremal Skolem sequences \\
%\hline
%1 & 1 \\
%2 & 1 \\
%3 & 3 \\
%4 & 7 \\
%5 & 23 \\
%6 & 83 \\
%7 & 405 \\
%8 & 2113 \\
%9 & 12657 \\
%10 & 82297 \\
%11 & 596483 \\
%12 & 4698655 \\
%13 & 40071743 \\
%\hline
%\end{tabular}
%%\end{center}
%\label{nrextskol}
%\end{table}

\begin{table}[!h]
%\begin{center}
\caption{Number of perfect extremal Skolem-type sequences.}
\begin{tabular*}{13cm}{|c@{\extracolsep{\fill}}|ccccccccc|}
\hline
Order & 1 & 2 & 3 & 4 & 5 & 6 & 7 & 8 & 9 \\
\hline
$\#$Sequences $\;$ & 1 & 1 & 3 & 7 & 23 & 83 & 405 & 2113 & 12675\\
\hline
\end{tabular*}
\begin{tabular*}{13cm}{|c@{\extracolsep{\fill}}|cccc|}
\hline
Order & 10 & 11 & 12 & 13  \\
\hline
$\#$Sequences $\;$ & 82297 & 596483 & 4698655 & 40071743 \\
\hline
\end{tabular*}
%\end{center}
\label{nrextskol}
\end{table}

The reader may have noticed that the number of perfect Skolem-type sequences of a given order seems to be odd. 
This is because perfect Skolem-type sequences are paired via the operation of reversing a sequence, except for the one sequence 
(of every order) of the form
$\dots 75311357 \dots$, which is invariant under reversal.

Moreover, the following conjecture is strongly suggested by the data in Table \ref{nrskol}.
\begin{conj}
The number of perfect Skolem-type sequences of order $n$ is always of the form $4k +1$.
\label{countconj}
\end{conj}

Another interesting enumerative question is that of the number of perfect Skolem sets of a given order. If a proof of Conjecture 
\ref{perfconj} were found it might be possible to find a closed formula for the number of perfect Skolem sets of any given order. 
Table \ref{nrskolset} contains the 
number of perfect Skolem sets of orders $1$ to $20$.

%\begin{table}[!h]
%%\begin{center}
%\caption{Number of perfect Skolem sets.}
%\label{nrskolset}
%\begin{tabular}{|c|l|}
%\hline
%Order & $\#$perfect Skolem sets \\
%\hline
%1 & 1 \\
%2 & 1 \\
%3 & 3 \\
%4 & 11 \\
%5 & 35 \\
%6 & 114 \\
%7 & 407 \\
%8 & 1486 \\
%9 & 5414 \\
%10 & 19923 \\
%11 & 74230 \\
%12 & 278462 \\
%13 & 1049318 \\
%\hline
%\end{tabular}
%%\end{center}
%\end{table}

\begin{table}[!h]
%\begin{center}
\caption{Number of perfect Skolem sets.}
\begin{tabular*}{13cm}{|c@{\extracolsep{\fill}}|ccccccccc|}
\hline
Order & 1 & 2 & 3 & 4 & 5 & 6 & 7 & 8 & 9 \\
\hline
$\#$Sets $\;$ & 1 & 1 & 3 & 11 & 35 & 114 & 407 & 1486 & 5414\\
\hline
\end{tabular*}
\begin{tabular*}{13cm}{|c@{\extracolsep{\fill}}|cccccc|}
\hline
Order & 10 & 11 & 12 & 13 & 14 & 15 \\
\hline
$\#$Sets $\;$ & 19923 & 74230 & 278462 & 1049318 & 3972395 & 15101658\\
\hline
\end{tabular*}
\begin{tabular*}{13cm}{|c@{\extracolsep{\fill}}|ccccc|}
\hline
Order & 16 & 17 & 18 & 19 & 20 \\
\hline
$\#$Sets $\;$ & 57607431 & 220391316 & 845366406 & 3250192681 & 12521965697 \\
\hline
\end{tabular*}
%\end{center}
\label{nrskolset}
\end{table}

%\begin{table}[!h]
%%\begin{center}
%\caption{Number of perfect extremal multi Skolem sets (A019589).}
%\label{nrextmultiskolset}
%\begin{tabular}{|c|l|}
%\hline
%Order & $\#$perfect extremal multi Skolem sets \\
%\hline
%1 & 1 \\
%2 & 2 \\
%3 & 5 \\
%4 & 16 \\
%5 & 59 \\
%6 & 246 \\
%7 & 1105 \\
%8 & 5270 \\
%9 & 26231 \\
%10 & 135036 \\
%11 & 713898 \\
%12 & 3857113 \\
%%13 & ? \\
%\hline
%\end{tabular}
%%\end{center}
%\end{table}

%\begin{table}[!h]
%%\begin{center}
%\caption{Number of perfect extremal Skolem sets (A000571).}
%\label{nrextskolset}
%\begin{tabular}{|c|l|}
%\hline
%Order & $\#$perfect extremal Skolem sets \\
%\hline
%1 & 1 \\
%2 & 1 \\
%3 & 2 \\
%4 & 4 \\
%5 & 9 \\
%6 & 22 \\
%7 & 59 \\
%8 & 169 \\
%9 & 490 \\
%10 & 1486 \\
%11 & 4639 \\
%12 & 14805 \\
%13 & 48107 \\
%\hline
%\end{tabular}
%%\end{center}
%\end{table}

\section{Restricted permutations}
\label{restrictedpermutations}
In this section we investigate connections between Skolem-type sequences and various problems on permutations. As we will see,
extremal Skolem-type sequences will be of particular interest with respect to these problems.

\subsection{Involutions and Skolem-type sequences}
\label{invskol}
We begin by stating some general observations on the connection between Skolem-type sequences and involutions.
An involution is a permutation which is its own inverse. A fixed point free involution is an involution where no element is mapped to 
itself. Another way to express this is that a fixed point free involution is a permutation such that when written in cycle notation 
it only has cycles of length $2$. 
It is clear that the proof of Theorem \ref{nrmulti} gives us a one-to-one correspondence between fixed point free 
involutions in $S_{2n}$ and perfect multi Skolem-type sequences of order $n$.

If we restrict the fixed point free involutions and require that all the distances between the transposed elements must be distinct, 
then we have a one-to-one correspondence between these restricted involutions in $S_{2n}$ and 
perfect non-multi Skolem-type sequences of order $n$. For example, the involution $(15)(23)(46)$ satisfies this condition because 
$|1-5| = 4,$ $|2-3| = 1$ and $|4-6| = 2$. Thus, Table \ref{nrskol} gives us the number of these restricted involutions in $S_{2n}$ 
when $n \leq 13$.

\begin{figure}[!h]
\begin{center}
\includegraphics[scale = 0.6]{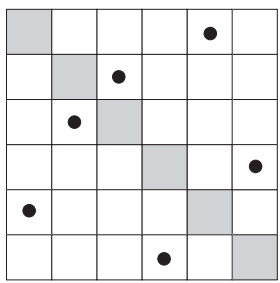}
\\
\end{center}
\caption{The permutation matrix corresponding to the fixed point free involution $(15)(23)(46)$ (perfect Skolem-type sequence $411242$). 
The shaded squares indicate forbidden positions.}
\label{fixfreeinv}
\end{figure}

The permutation matrix of a fixed point free involution is symmetric and has the northwest to southeast main diagonal empty. 
The restriction on the involutions has a nice interpretation in terms of permutation matrices, namely that every northwest 
to southeast diagonal of the matrix must have at most one entry. Two elements on the same northwest to southeast diagonal imply 
that two of the transpositions in the involution are of the form $(ab)$ and $((a+k)(b+k))$, and of course $|a-b| = |(a+k)-(b+k)|$, 
violating the condition. See Figure \ref{fixfreeinv} for an example of a permutation matrix representing an involution satisfying
the condition.

\begin{figure}[!h]
\begin{center}
\includegraphics[scale = 0.6]{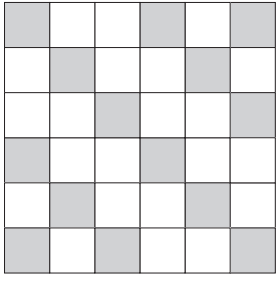}
\\
\end{center}
\caption{A matrix such that elements can be added to it to yield a symmetric permutation matrix where every northwest to 
southeast diagonal contains at most one element if and only if the set $A = \{1,2,4\}$ is a Skolem set. The shaded squares 
indicate forbidden positions.}
\label{perffixfreeinv}
\end{figure}

Given a set $A$ of cardinality $n$ we construct a $2n \times 2n$ matrix $M$ which have certain forbidden positions, namely the 
main northwest to southeast diagonal and those $m_{ij}$ where $|i-j| \notin A$. Elements can then be added to $M$ to yield a
 symmetric permutation matrix with at most one element on every northwest to southeast diagonal if and only if $A$ 
is a Skolem set. See Figure \ref{perffixfreeinv} for an example.

\subsection{Permutations and extremal Skolem-type sequences}
\label{permextremalskol}
%Permutations are among the richest and most fundamental objects in discrete mathematics. Many different aspects of permutations has
%been investigated, far to many to give account for here. From a combinatorial perspective much of the research has been concentrated 
%on enumerative questions. 
In this section we will consider an interesting existence question for permutations and its relation to perfect 
extremal Skolem-type sequences. 

To the best of our knowledge the following type of fundamental existence question for permutations has not been studied before.
Given a sequence of say $5$ elements, is there a permutation of the elements such that one element is moved $4$ steps to the right, one
element is moved $1$ step to the right, one element remains in its original position, one element is moved $2$ step to the left, 
and one element is moved $3$ steps to the left? 

We define the problem formally as follows.
\begin{defn}
To every permutation $\pi$ in $S_n$ we associate a displacement pattern $\alpha_\pi$ defined as follows, 
let $\alpha_{\pi}$ be the multiset of differences $\{ \pi(i) - i \; | \; 1 \leq i \leq n\}$ sorted in decreasing order.
Given a displacement pattern $\alpha$ we want to determine
whether there exists a permutation $\pi$ having the displacement pattern $\alpha_\pi = \alpha$.  
\label{ppatterndef}
\end{defn}

Next we present the connection between the problem of permutation displacements, as defined above, and Skolem-type sequences.
\begin{thm}
Given a displacement pattern $\alpha = (a_1,a_2,\dots,a_n)$,
then there exists a permutation $\pi$ in $S_n$ having 
the displacement pattern $\alpha_\pi = \alpha$ if and only if the multiset $A = \{a_1 + n, a_2 + n, \dots, a_n + n\}$ is a perfect
extremal multi Skolem set.
\label{permcorresp}
\end{thm}
\begin{pf}
First of all note that the proof of Theorem \ref{nrmaxmulti} gives us a one-to-one correspondence between perfect extremal 
multi Skolem-type sequences of order $n$ and permutations in $S_n$.

We begin by proving the if part.
If the multiset $A = \{a_1 + n, a_2 + n, \dots, a_n + n\}$ is a perfect
extremal multi Skolem set, then we know that the set $\{1,2,\dots, 2n\}$ can be partitioned into the differences in $A$. Moreover, since
$A$ is extremal we know that this partition $a_i + n = s_i - t_i$, ($i = 1,2,\dots,n$) has the property that 
$t_i \in \{1,2,\dots,n\}$ and $s_i \in \{n+1,n+2,\dots,2n\}$, ($i = 1,2,\dots,n$). Now consider the permutation $\pi$ on 
$\{1,2,\dots,n\}$ defined as 
$\pi(t_i) = s_i - n$, ($i = 1,2,\dots,n$). The corresponding displacement pattern is 
$\alpha_\pi = (s_1 - n - t_1, s_2 - n - t_2, \dots, s_n - n - t_n) = (a_1, a_2, \dots, a_n)$.

Now to the only if part.
If $\pi$ is a permutation having displacement pattern $\alpha_\pi = (a_1,a_2,\dots,a_n)$, then
$(a_1,a_2,\dots,a_n) = (\pi(t_1) - t_1, \pi(t_2) - t_2, \dots, \pi(t_n) - t_n)$ where 
$\{t_1, t_2, \dots, t_n\} = \{1,2,\dots,n\}$. Hence, $(\pi(t_i) + n) - t_i$, ($i = 1,2,\dots,n$) is a partition of 
$\{1,2,\dots, 2n\}$ into the differences in $\{a_1 + n, a_2 + n, \dots, a_n + n\}$. Thus, 
$A = \{a_1 + n, a_2 + n, \dots, a_n + n\}$ is a perfect
extremal multi Skolem set.
\qed
\end{pf}
Note that by the preceding theorem, the existence question for permutation displacements in Definition \ref{ppatterndef} is merely a
special case of the existence question for perfect multi Skolem-type sequences.

The preceding theorem together with the necessary conditions for the existence of perfect multi Skolem-type sequences in 
Corollary \ref{necperf} imply the following necessary conditions for the existence of permutations satisfying a given
displacement pattern.
\begin{cor}
Given a displacement pattern $\alpha = (a_1,a_2,\dots,a_n)$, where $a_1 \geq a_2 \geq \dots \geq a_n$. If there exists a permutation $\pi$ in $S_n$ having displacement pattern
$\alpha_\pi = \alpha$, then $\sum^m_{i=1} a_i \leq m(n-m)$ for each 
$1 \leq m \leq n$, with equality holding when $m=n$ (the density condition). Further more, if $n$ is even, then the number of even 
elements in $\alpha$ is even, and if $n$ is odd, then 
the number of odd elements in $\alpha$ is even (the parity condition).
\label{permnec}
\end{cor}
These necessary conditions are not sufficient, the displacement pattern $\alpha = (2,2,2,-2,-2,-2)$ is a minimal counterexample.

From Theorem \ref{permcorresp} we make the following obvious but important observation.
\begin{cor}
Given a displacement pattern $\alpha = (a_1,a_2,\dots,a_n)$ where all $a_i$'s are distinct, i.e., $a_1 > a_2 > \dots > a_n$,
then there exists a permutation $\pi$ in $S_n$ having 
the displacement pattern $\alpha_\pi = \alpha$ if and only if the set $A = \{a_1 + n, a_2 + n, \dots, a_n + n\}$ is a perfect
extremal Skolem set.
\label{permcorrespnonmulti}
\end{cor}

Hence, the following conjecture is just a special case of Conjecture \ref{perfconj}.
\begin{conj}
Given a displacement pattern $\alpha = (a_1,a_2,\dots,a_n)$ where all $a_i$'s are distinct, then the necessary conditions in
Theorem \ref{permnec} are also sufficient for the existence of a permutation $\pi$ having displacement pattern $\alpha_{\pi} = \alpha$.
\label{permconj}
\end{conj}
%The conjecture has been verified by computer search up to $n=13$.

By the group structure of permutations we get the following result.
\begin{thm}
Given a displacement pattern $\alpha = (a_1,a_2,\dots,a_n)$,
then there exists a permutation $\pi$ in $S_n$ having 
the displacement pattern $\alpha_\pi = \alpha$ if and only if there exists a permutation $\rho$ in $S_n$ having the displacement
pattern $\alpha_\rho = (-a_n,-a_{n-1},\dots,-a_1)$.
\label{symmetryperm}
\end{thm}
\begin{pf}
Given a permutation $\pi$ in $S_n$ having displacement pattern $\alpha_\pi = (a_1,a_2,\dots,a_n)$, then there is a (unique) permutation
$\pi^{-1}$ such that $\pi \circ \pi^{-1} = \iota$, where $\iota$ is the identity permutation. Since $\iota$ has the displacement pattern
$(0,0,\dots,0)$ we get that $\pi^{-1}$ has the displacement pattern $(-a_n,-a_{n-1},\dots,-a_1)$. To prove the other direction, just
reverse the roles of $\pi$ and $\pi^{-1}$.
\qed
\end{pf}

The preceding theorem together with Theorem \ref{permcorresp} immediately gives us the following result.
\begin{cor}
$A = \{a_1, a_2, \dots, a_n\}$ is a perfect
extremal multi Skolem set if and only if $A = \{2n-a_n, 2n-a_{n-1}, \dots, 2n-a_1\}$ is a perfect extremal multi Skolem set.
\label{symmetryseq}
\end{cor}

Most of the previous research on permutation displacements have focused on enumerative aspects. For example, Lehmer \cite{Lehmer}
studied the following problems. Determine the number of permutations $\pi$ in $S_n$ such that
\begin{enumerate}
\item
no element is moved more than $k$ positions left or right, that is
$|\pi(i) - i| \leq k$ $(i = 1,\dots,n)$, denoted $P_1^k$; 
\item
the last element is moved to the first position, all other elements are moved right not more than $k$ positions, that is
$\pi(n)-n = -n+1$ and $1 \leq \pi(i) -i \leq k$ $(i = 1,\dots,n-1)$, denoted $P_2^k$;
\item
no element is moved more than $k$ positions left or right but each element must move, that is
$1 \leq |\pi(i) - i| \leq k$ $(i = 1,\dots,n)$, denoted $P_3^k$.
\end{enumerate}
Lehmer gives solutions to $P_1^1,P_1^2,P_1^3,P_2^k,P_3^1,$ and $P_3^2$ in terms of generating functions, e.g.,
$P_2^k$ is given by the generating function $\frac{1-x-x^2+x^{k+1}}{1-2x+x^k}$. Note that $P_3^n$ is the well known 
derangement problem which has the solution $n! \sum_{i=0}^n \frac{(-1)^i}{i!}$ \cite{stan1}.

It is not hard to see that the correspondence between displacement patterns and perfect extremal multi Skolem sets in 
Theorem \ref{permcorresp} is solution preserving in the sense that the number of permutations having the displacement pattern 
$\alpha = (a_1,a_2,\dots,a_n)$ equals the number of perfect extremal multi Skolem-type sequences that can be generated from
the multiset $A = \{a_1 + n, a_2 + n, \dots, a_n + n\}$. Hence, enumerative results on permutation displacements can be transferred to
enumerative results on perfect extremal multi Skolem-type sequences, and vice versa. For example, the problem of counting the 
number of perfect 
extremal multi Skolem-type sequences of order $n$ having no occurrence of $n$ is just a reformulation of the derangement problem, 
and thus the solution is $n! \sum_{i=0}^n \frac{(-1)^i}{i!}$. Moreover, our count of the number of perfect extremal (non-multi) 
Skolem-type sequences of order $n$ ($n \leq 13$) in Table \ref{nrextskol} gives us the number of permutations in $S_n$ ($n \leq 13$) 
where no two elements are moved the same number
of positions in the same direction, or equivalently the number of $n \times n$ permutation matrices ($n \leq 13$) where each
northwest to southeast diagonal contains at most one element. The problem of finding a closed formula (or even a recurrence relation) for
the number of these sequences/permutations seems like a very interesting and challenging problem.

We conclude this section with an enumerative result that is an easy consequence of the proof of Theorem \ref{symmetryperm} and 
Corollary \ref{symmetryseq}.
\begin{cor}
The number of perfect extremal multi Skolem-type sequences that can be generated from $A =  \{a_1, a_2, \dots, a_n\}$ is equal to the number
of perfect extremal multi Skolem-type sequences that can be generated from $A = \{2n-a_n, 2n-a_{n-1}, \dots, 2n-a_1\}$.
\label{nrsymmetryseq}
\end{cor}

\section{Constructions}
\label{sectionconstr}
In this section we show that every perfect multi Skolem set of order $n$ gives rise to a new perfect extremal multi Skolem set
of order $2n$. We also show that every $k$-extended multi Skolem set of order $n$ gives rise to a new perfect extremal multi Skolem
set of order $2n+1$. We show how these results can be used to automatically transform many existence results for perfect Skolem sets and
$k$-extended Skolem sets into existence results for perfect extremal Skolem sets. Moreover, we give a complete solution to
the existence question for extremal near-Langford sequences.

\begin{thm}
If $A = \{ a_1, a_2, \dots, a_n\}$ is a perfect multi Skolem set, then 
$B = \{2n + a_1, 2n + a_2, \dots, 2n + a_n, 2n - a_n, 2n - a_{n-1}, \dots, 2n - a_1\}$ is a perfect extremal multi Skolem set.
\label{extfromskol}
\end{thm}
\begin{pf}
If $A = \{ a_1, a_2, \dots, a_n\}$ is a perfect multi Skolem set, then we know that there exists a partition 
of $\{1,2,\dots,2n\}$ into the differences in $A$. As we have seen in Section \ref{invskol} this partition can be seen as a 
fixed point free involution in $S_{2n}$. The displacement pattern corresponding to this involution is 
$(a_1, a_2, \dots, a_n, -a_n, -a_{n-1}, \dots, -a_1)$. Now by Theorem \ref{permcorresp}, we get that 
$B = \{2n + a_1, 2n + a_2, \dots, 2n + a_n, 2n - a_n, 2n - a_{n-1}, \dots, 2n - a_1\}$ is a perfect extremal multi Skolem set.
\qed
\end{pf}
Note that if $A$ is a perfect (non-multi) Skolem set, then $B$ is a perfect extremal (non-multi) Skolem set.

\begin{thm}
If $A = \{ a_1, a_2, \dots, a_n\}$ is a $k$-extended multi Skolem set, then 
$B = \{2n + 1 + a_1, 2n + 1 + a_2, \dots, 2n + 1 + a_n, 2n + 1, 2n + 1 - a_n, 2n + 1 - a_{n-1}, \dots, 2n + 1 - a_1\}$ is a perfect extremal multi Skolem set.
\label{extfromskol2}
\end{thm}
\begin{pf}
If $A = \{ a_1, a_2, \dots, a_n\}$ is a $k$-extended multi Skolem set, then we know that there exists a partition 
of $\{1,2,\dots,2n,2n+1\} \setminus k$ into the differences in $A$. As we have seen in Section \ref{invskol} this 
partition can be seen as an involution in $S_{2n+1}$ having exactly one fixed point. The displacement pattern 
corresponding to this involution is 
$(a_1, a_2, \dots, a_n,0, -a_n, -a_{n-1}, \dots, -a_1)$. Now by Theorem \ref{permcorresp}, we get that 
$B = \{2n + 1 + a_1, 2n + 1 + a_2, \dots, 2n + 1 + a_n, 2n + 1, 2n + 1 - a_n, 2n + 1 - a_{n-1}, \dots, 2n + 1 - a_1\}$ 
is a perfect extremal multi Skolem set.
\qed
\end{pf}
Note that if $A$ is a $k$-extended (non-multi) Skolem set, then $B$ is a perfect extremal (non-multi) Skolem set.

Any perfect multi Skolem set is trivially also a $0$-extended multi Skolem set. Hence, the following corollary follows directly 
from Theorem \ref{extfromskol2}.
\begin{cor}
If $A = \{ a_1, a_2, \dots, a_n\}$ is a perfect multi Skolem set, then 
$B = \{2n + 1 + a_1, 2n + 1 + a_2, \dots, 2n + 1 + a_n, 2n + 1, 2n + 1 - a_n, 2n + 1 - a_{n-1}, \dots, 2n + 1 - a_1\}$ is a perfect extremal multi Skolem set.
\label{extfromskol3}
\end{cor}

The preceding theorems together with known existence results for Skolem-type sequences give us several new existence results for
extremal Skolem-type sequences essentially for free. For example, the following result follows from the existence results for
near-Skolem sequences and hooked near-Skolem sequences in \cite{Shal}.
\begin{cor}
$A= \{3n-1,3n-2, \dots,n-1\} \setminus \{2n-1+m, 2n-1-m\}$ is a perfect extremal Skolem set for all $1 \leq m < n$.
\end{cor}
\begin{pf}
Shalaby proved in \cite{Shal} that a set of the form $A = \{1,2,\dots,m-1,m+1,\dots,n\}$, where $m$ is even, is a perfect 
Skolem set for $n \equiv 2,3 \pmod{4}$, and a ($2n-3$)-extended Skolem set for $n \equiv 0,1 \pmod{4}$.
Hence, by Theorem \ref{extfromskol2} and Corollary \ref{extfromskol3}, we get that 
$A= \{3n-1,3n-2, \dots,n-1\} \setminus \{2n-1+m, 2n-1-m\}$, where $m$ is even, is a perfect extremal Skolem set.
Shalaby also proved that a set of the form $A = \{1,2,\dots,m-1,m+1,\dots,n\}$, where $m$ is odd, is a perfect Skolem set for
$n \equiv 0,1 \pmod{4}$, and a ($2n-3$)-extended Skolem set for $n \equiv 2,3 \pmod{4}$. Again, by Theorem \ref{extfromskol2} and 
Corollary \ref{extfromskol3}, we get that $A= \{3n-1,3n-2, \dots,n-1\} \setminus \{2n-1+m, 2n-1-m\}$, where $m$ is odd, is 
a perfect extremal Skolem set.
\qed
\end{pf}

We give some more examples of the same flavour.
\begin{cor}
$A$ is a perfect extremal (multi) Skolem set if
\begin{enumerate}
\item
$A= \{3n-2, 3n-3,\dots,n-2 \} \setminus \{2n-2+m, 2n-2, 2n-2-m\}$ 
and $n \equiv 0,1 \pmod{4}$ and $m$ is odd, 
or $n \equiv 2,3 \pmod{4}$ and $m$ is even; or
\item
%$A=\{(3n)^m,(3n-1)^m,\dots,(2n+2)^m,(2n+1)^m,(2n-1)^m,(2n-2)^m,\dots,(n+1)^m,n^m\}$ 
$A= \{ (2nm+n)^m, \dots, (2nm+1)^m,  (2nm-1)^m, \dots, (2nm-n)^m \}$ and $n \equiv 0,1 \pmod{4}$, or
$n \equiv 2,3 \pmod{4}$ and $m$ is even.
\end{enumerate}
\end{cor}
\begin{pf}
The proof of (1) follows directly from Theorems \ref{nearskolem} and \ref{extfromskol}, and the proof of
(2) follows directly from Theorems \ref{mfold} and \ref{extfromskol}.
\qed
\end{pf}

Next we show how pairs of perfect extremal multi Skolem sets give rise to new perfect extremal multi Skolem sets.  
\begin{thm}
If $A = \{ a_1, a_2, \dots, a_n\}$ and $B = \{ b_1, b_2, \dots, b_m\}$ are perfect extremal multi Skolem sets, then so are
\begin{enumerate}
\item
$C= A \cup \{b_1+2n, b_2 + 2n, \dots, b_m + 2n\}$,
\item
$D= B \cup \{a_1+2m, a_2 + 2m, \dots, a_n + 2m\}$, and
\item
$E= \{a_1+m, a_2 + m, \dots, a_n + m\} \cup \{b_1 + n, b_2 + n, \dots, b_m + n\}$.
\end{enumerate}
\label{extremalrecursive}
\end{thm}
\begin{pf}
If $A$ and $B$ are perfect extremal multi Skolem sets, then they can be used to construct the perfect extremal multi Skolem-type sequences
$S_A$ and $S_B$ respectively. If we split $S_A$ in the middle we obtain two sequences $S_A^L$ and $S_A^R$ consisting of the $n$ leftmost and rightmost elements in $S_A$ respectively. 
It is crucial to realize that since we are working with extremal sequences we know that any two elements that are paired by the partition
of $\{1,\dots,2n\}$ into two-tuples induced by $S_A$ will have one element in $S_A^L$ and one element in $S_A^R$.
We define $S_B^L$ and $S_B^R$ analogously. 
Let $S + k$, where $S$ is a sequence and $k$ a natural number, denote the sequence obtained by increasing each element in $S$ by $k$.

Now, to prove that $C= A \cup \{b_1+2n, b_2 + 2n, \dots, b_m + 2n\}$ is a perfect extremal multi Skolem set, we give a perfect extremal multi Skolem-type sequence $S_C$ which
can be generated from $A \cup \{b_1+2n, b_2 + 2n, \dots, b_m + 2n\}$. Let $S_C = (S_B^L+2n)S_A(S_B^R+2n)$. Similarly 
$S_D = (S_A^L+2m)S_B(S_A^R+2m)$ is a perfect extremal multi Skolem-type sequence that can be generated from $B \cup \{a_1+2m, a_2 + 2m, \dots, a_n + 2m\}$. Finally, $S_E = (S_A^L + m)(S_B^L+n)(S_A^R+m)(S_B^R+n)$ is a perfect extremal multi Skolem-type sequence which can be created from
$\{a_1+m, a_2 + m, \dots, a_n + m\} \cup \{b_1 + n, b_2 + n, \dots, b_m + n\}$.
\qed
\end{pf}
We believe that the preceding theorem together with the large number of existence results for perfect extremal Skolem-type sequences, that can automatically be obtained from existence results for (not necessarily extremal) Skolem-type sequences, make them particularly well suited as building blocks for proving existence results for perfect (not necessarily extremal) Skolem-type sequences.

%\subsection{Near Langford sequences with small defects}
%blablabla k extended langford with small defects blabla bla definition of near langford with small defects
The last result of this section is a complete solution to the existence question for extremal near-Langford sequences.
Recall that a Langford sequence is perfect Skolem-type sequence that can be generated from a set of the form 
$A = \{a,a+1,\dots, b\}$, where $a$ is usually referred to as the {\em defect\/} of the sequence/set.
Skolem solved the existence question for Langford sequences of defect $1$ in \cite{Skol}; Priday solved it for defect $2$ \cite{Priday};
Bermond et al. solved it for defect $3$ and $4$ \cite{BBG76}; before Simpson gave a complete solution for all defects \cite{Simp}.
A near-Langford sequence of defect $a$ is perfect Skolem-type sequence that can be generated from a set of the form
$A=\{a,a+1,\dots,b\} \setminus \{m\}$, where $a<m<b$. Shalaby solved the existence question for near-Langford sequences of defect $1$ in
\cite{Shal} (see also Theorem \ref{nearskolem}), the problem is open for all other defects. 
\begin{thm}
A set $A=\{a,a+1,\dots,b\} \setminus \{m\}$, where $a<m<b$, is a perfect extremal Skolem set if and only if $b = 3a$ and $m = (a + b)/2$.
\label{extnearlang}
\end{thm}
\begin{pf}
We begin by proving the only if part. Assume (with the aim of reaching a contradiction) that $A=\{a,a+1,\dots,b\} \setminus \{m\}$ is a perfect extremal Skolem set, but $b \neq 3a$. Two cases emerge
\begin{enumerate}
\item
$b > 3a$, i.e., $b = 3a + k$ where $0<k$, or 
\item
$b < 3a$, i.e., $b = 3a - k$ where $0<k<2a$.
\end{enumerate}
Recall from the definition of perfect extremal Skolem sets that any perfect extremal Skolem set $A=\{a_1,a_2,\dots,a_n\}$ must satisfy
$$\sum^n_{i=1} a_i = n^2.$$
Hence, we have that $A=\{a,a+1,\dots,b\} \setminus \{m\}$ must satisfy
\begin{equation}
(\sum_{i=a}^{b} i) - m = (b-a)^2.
\label{equationsum}
\end{equation}

Now, consider the first case above where $b = 3a + k$, and $0<k$. Writing the sum in closed form and substituting $b$ by $3a +k$ in
(\ref{equationsum}) above gives us
\begin{equation}
\frac{(4a+k)(2a+k+1)}{2}-m = (2a+k)^2.
\label{equationintermediate}
\end{equation}
By straightforward formula manipulations, we get
\begin{equation}
4a-2ka+k-k^2 = 2m.
\label{equationfinal}
\end{equation}
The lhs of (\ref{equationfinal}) above is less than or equal to $2a$ for all $k>0$. This implies that $m \leq a$ for all $k>0$, 
which is a contradiction with the fact that $a<m$.

Now, consider the second case above where $b = 3a - k$, and $0<k<2a$. 
By arguments analogous to those in the previous case, we deduce that the equality
\begin{equation}
4a + 2ka -k -k^2 = 2m
\label{equationcase2}
\end{equation}
must hold.
We know that $m$ must be less than $b$ and since $b = 3a-k$ we can deduce from \ref{equationcase2} that the inequality 
\begin{equation}
2(3a -k) > 4a + 2ka -k -k^2
\label{final2}
\end{equation}
must hold. By simplifying (\ref{final2}), we get $2a > 2ak - k^2 + k$, and by elementary calculus it is easy to see that $2ak -k^2 + k \geq 2a$ for all $k$ in the interval $0<k<2a$. Hence, again we have a contradiction and it follows that $b=3a$.

If $b=3a$, then we substitute $b$ by $3a$ in (\ref{equationsum}). A simple calculation verifies that $m = 2a = (a+b)/2$, which concludes
the proof of the only if part of the theorem. To prove the if part of the theorem, we show in Table \ref{extrnearl} how to construct
perfect extremal Skolem-type sequences from all sets of the form $A=\{a,a+1,\dots,b\} \setminus \{m\}$, where $a<m<b$, $b = 3a$, and $m = (a + b)/2$.

\begin{table}[!h]
\begin{center}
%\caption{Perfect extremal Skolem sequences generated from $A=\{a,a+1,\dots,b\} \setminus \{m\}$, where $a<m<b$, $b = 3a$, and $m = (a + b)/2$}
\caption{}
\begin{tabular*}{13.8cm}{c@{\extracolsep{\fill}}cccc}
%\begin{tabular}{ccccc}
\hline
 & $a_i$ & $b_i$ & $|b_i-a_i|$ & $j\in$ \\
\hline
(1) & $2a-2j$ & $3a-j$ & $a+j$ & $[0,a-1]$ \\
(2) & $4a+1-2j$ & $5a+1-j$ & $a+j$ & $[a+1,2a]$  \\
\hline
\end{tabular*}
\end{center}
\label{extrnearl}
\end{table}
\qed
\end{pf}

%\section{Complexity}
\section{Final remarks}
After having completed this work, the author learned about the recent results in \cite{Postnikov,Brunetti}. Both these papers
study a problem that is easily seen to be equivalent to the existence question for perfect extremal multi Skolem-type sequences.
Several interesting connections between this problem and other combinatorial objects are presented in \cite{Postnikov}. Moreover,
Conjecture 4.2 in \cite{Postnikov} is easily seen to be a reformulation of a special case of our Conjecture \ref{perfconj} (in fact,
it is equivalent to our Conjecture \ref{extconj}).

The authors of \cite{Brunetti} study the computational complexity of deciding whether a given set $A$ is a perfect extremal 
multi Skolem set (or rather a question equivalent to this one). 
They manage to prove that this problem is NP-complete, a result which clearly
subsumes the NP-completeness result for generalized Skolem sequences in \cite{Nordh}. Note that since the existence question
for perfect extremal multi Skolem-type sequences is just a special case of the existence question for perfect multi Skolem-type
sequences, it follows that the problem of deciding whether a set $A$ is a perfect multi Skolem set is also NP-complete.
The NP-completeness result in \cite{Brunetti} is important
since it offers a very good explanation of why we have not been able to extend Conjecture \ref{perfconj} to the case where
$A$ is a multiset. More specifically, the results in \cite{Brunetti} imply that there can be no simple necessary and 
sufficient conditions for the existence of perfect multi Skolem-type sequences unless P$=$NP (where simple means checkable in polynomial time).

\begin{ack} %{Acknowledgements}
This work is in part a result of the author's master's thesis \cite{Nordh}. The author thanks the supervisor Svante Linusson for 
proposing this interesting problem; Daniel Hedberg, Peter Jonsson, Václav Linek, and Nabil Shalaby for valuable remarks on earlier 
versions of the manuscript; the anonymous referees for their comments that have improved the exposition in many ways.
\end{ack}

%\begin{thebibliography}{00}

% \bibitem{label}
% Text of bibliographic item

% notes:
% \bibitem{label} \note

% subbibitems:
% \begin{subbibitems}{label}
% \bibitem{label1}
% \bibitem{label2}
% If there is a note, it should come last:
% \bibitem{label3} \note
% \end{subbibitems}

%\bibitem{}

%\end{thebibliography}

\bibliography{references}
\bibliographystyle{plain}
\end{document}